\documentclass[12pt]{amsart}
\makeatletter

\newtheorem{thm}{Theorem}[section]
\def\statetheorem{\@ifnextchar[{\@statetheorem}{\nr@statetheorem}}
\long\def\@statetheorem[#1]#2{\begin{thm}\label{#1}#2\end{thm}}
\long\def\nr@statetheorem#1{\begin{thm}#1\end{thm}}
\def\statetheorempf{\@ifnextchar[{\@statetheorempf}{\nr@statetheorempf}}
\long\def\@statetheorempf[#1]#2{\begin{thm}\label{#1}#2\end{thm}\proof}
\long\def\nr@statetheorempf#1{\begin{thm}#1\end{thm}\proof}

\newtheorem{lmma}{Lemma}[section]
\def\statelemma{\@ifnextchar[{\@statelemma}{\nr@statelemma}}
\long\def\@statelemma[#1]#2{\begin{lmma}\label{#1}#2\end{lmma}}
\long\def\nr@statelemma#1{\begin{lmma}#1\end{lmma}}
\def\statelemmapf{\@ifnextchar[{\@statelemmapf}{\nr@statelemmapf}}
\long\def\@statelemmapf[#1]#2{\begin{lmma}\label{#1}#2\end{lmma}\proof}
\long\def\nr@statelemmapf#1{\begin{lmma}#1\end{lmma}\proof}

\newtheorem{crlry}{Corollary}[section]
\def\statecorollary{\@ifnextchar[{\@statecorollary}{\nr@statecorollary}}
\long\def\@statecorollary[#1]#2{\begin{crlry}\label{#1}#2\end{crlry}}
\long\def\nr@statecorollary#1{\begin{crlry}#1\end{crlry}}
\def\statecorollarypf{\@ifnextchar[{\@statecorollarypf}{\nr@statecorollarypf}}
\long\def\@statecorollarypf[#1]#2{\begin{crlry}\label{#1}#2\end{crlry}\proof}
\long\def\nr@statecorollarypf#1{\begin{crlry}#1\end{crlry}\proof}

\def\note{\@ifnextchar[{\@note}{\@note[Note]}}
\def\@note[#1]{\par\medskip\noindent{\textbf{#1:} }}
\def\C{\mathbb{C}}

\def\N{\mathbb{N}}
\def\bibref[#1]{\cite{#1}}

\let\real@bibitem\bibitem
\def\bibitem[#1]{\real@bibitem{#1}}

\makeatother

\def\spec{\textup{Spec}\,}
\def\A{\mathcal{A}}
\def\Dd{{\mathcal{D}}}
\def\PD{\Psi\Dd}
\def\R{\mathcal{R}_{\lambda}}

\def\v{\vec{\hbox{\textbf{v}}}}

\begin{document}

\title{Commutative Partial Differential Operators}

\author{Alex Kasman}
\address{\hskip-\parindent
        Department of Mathematics\\
College of Charleston\\
66 George Street\\
Charleston, SC 29424-0001}
\email{kasman@math.cofc.edu}

\author{Emma Previato}
\address{\hskip-\parindent
        Department of Mathematics\\
        Boston University\\
        111 Cummington Street\\
        Boston, MA 02215}

\email{ep@math.bu.edu}

\thanks{Research at MSRI is supported in part by NSF grant DMS-9701755.
The second author was a member of the Symbolic Calculation program and
wishes to express tremendous gratitude to organizer Michael Singer.
The first author was a member of the Random Matrices program and
wishes to thank all of the organizers of that program and also thank
Michael Singer for his assistance and guidance.}

\begin{abstract}
In one variable, there exists a satisfactory classification of
commutative rings of differential operators.  In
several variables, even the simplest generalizations seem to be
unknown and in this report we give examples and pose questions that
may suggest a theory to be developed.   In particular, we address the
existence of a ``spectral variety'' generalizing the spectral curve of
the one dimensional theory and the role of the differential resultant.
\end{abstract}

\maketitle

\section{Introduction}

In one variable, there exists a satisfactory classification of
rings of differential operators that are (maximal) commutative.  In
several variables, even the simplest generalizations seem to be
unknown and in this report we give examples and pose questions that
may suggest a theory to be developed. 
To motivate what we do, we briefly recall the 1-variable case and
selected results in several variables.

In the 1-variables case with analytic coefficients the classification
was found by Burchnall and Chaundy in the 1920s \bibref[BC1,BC2] by
essentially formal methods of differential algebra and some algebraic
function theory on the ``spectral curve'' $\spec \A$, where $\A$ is
the commutative ring, completed by a point ``at infinity''.  The
isospectral rings, roughlty speaking, form the Jacobi variety of this
projective algebraic curve and the ``Krichever map'' \bibref[Kr]
solves the inverse spectral problem explicitly.  For a ring $\A$
generated by a pair $L=\partial^n+u_{n-2}\partial^{n-2}+\cdots+u_0(x)$
(where $\partial=\partial/\partial x$) and $B=\partial^m+\cdots$ the
determinant of the $(n+m)\times(n+m)$
resultant matrix of $L-\lambda$ and $B-\mu$ \bibref[PrevSpec] is a
non-zero, polynomial $p\in\C[\lambda,\mu]$ such that $p(L,B)=0$.
Moreover, in the  ``rank 1 case'' $\textup{gcd}(n,m)=1$
\bibref[PrevSpec] one has that $p(\lambda,\mu)=0$ \textit{is} the
affine equation of the spectral curve.  A divisor of points on the
curve corresponding to the given element of the Jacobi variety is the
set of poles of the (normalized) gcd of $L-\lambda$ and $B-\mu$ at a
reference point $x=x_0$; the flow in $x$ is linear on the Jacobian.

In view of this, it is natural to aks at least the following questions
in several variables:
\begin{enumerate}
\item Is $\spec\A$ an affine variety of dimension $N$ for any maximal
commutative ring $\A$ of partial differential operators in $N$
variables?
\item For a ring $\A$ with generators $L_i$ ($1\leq i \leq N+1$), what
is the relationship between the differential resultant of
$L_i-\lambda_i$ and the equation of such a variety?
\end{enumerate}

More ambitiously, of course, one would ask for isospectral flows, a
compete classification, and the inverse spectral problem, in
increasing order of magnitude.  A beautiful generalization of
Burchnall-Chaundy theory was given by Nakayashiki \bibref[Nak1,Nak2]
using the Fourier-Mukai transform.  He associates commutative rings in
$N$ variables to a suitable $N$-dimensional abelian variety an some
additional choices, one for each element of its Picard variety.  But
these are not scalar operators, rather they have $(N!\times N!$)
matrix coefficients.

In this report, we give a negative answer to question (1) by using
techniques developed in \bibref[jmaa] and we offer some observations,
natural conjectures, and a strategy to treat (2).  We deal with the
scalar case only.  We include as an appendix the \textit{Mathematica}
code that can be used to compute the differential resultant of any set of
$N+1$ partial differential operators in $N$ variables, which is a
handy tool for checking properties on the available class of examples
such as \bibref[BK,CV].

\section{Geometric Structure of Maximal Commutative Rings}

It is well known that the commutative rings of ordinary differential
operators are finitely generated rings coordinate rings for algebraic
curves.  This forms the foundation of the Burchnall-Chaundy theory of
such rings \bibref[BC1,BC2,PrevSpec].  In contrast, very little is known
about the algebro-geometric structure of commutative rings of partial
differential operators.  A question of interest is to address the
problem of whether every commutative subring is contained in a
(larger) commutative ring requiring only a finite number of generators
over $\C$.  Such a result is relevant, for instance, to the
algebro-geometric investigations of quantum integrable systems
\bibref[BEG,HK]. 

In this section we use techniques from soliton theory
(namely Darboux transformation and Baker-Akhiezer functions) to study
the structure of \textit{certain} commutative rings of differential
operators.  
We are able to show that these rings are \textit{maximal} in the
sense that they are not contained in any larger commutative subrings
of the ring of differential operators.  This in itself is a difficult
task which is rarely achieved.  Then, in one example, we study the
structure of this ring more closely and note that it cannot be
constructed with only a finite number of generators over $\C$ and
hence is not the coordinate ring of an affine algebraic variety.

\note[Notation] 
Let $\Dd=\C(x_1,\cdots,x_n)[\partial_1,\ldots,\partial_n]$ be the ring
of rational coefficient differential operators in $n$ variables.
It will be useful to be able to refer also to
$\Dd_{0}=\C[\partial_1,\ldots,\partial_n]\subset \Dd$
 (the constant
coefficient differential operators) as well as
$\PD$, a ring of microdifferential operators containing $\Dd$ as well
as the inverse of the particular operator $K\in\Dd$ which will be
important below 
\bibref[Kashiwara] and $\PD_0$, the contant coefficient
microdifferential operators.

\subsection{What Commutes with Many Constant Coefficient Operators?}

In general, it is difficult to address the question of whether a
commutative ring of partial differential operators is maximal.  The
key which allows us to do it here is the observation that although the
centralizer of a single constant coefficient operator
$p(\partial_1,\ldots,\partial_n)$ will contain
non-constant coefficient operators, an operator commutes
all constant coefficient multiples of $p(\partial_1,\ldots,\partial_n)$ if and
only if it also has constant coefficients.

\begin{lemma}[lem:only-const]{Let $p\in\Dd_0$ be a non-zero constant coefficient
differential operator.  Then any operator $L\in\PD$ which commutes
with $p$ as well as all operators $P_i:=\partial_i\circ
p(\partial_1,\ldots,\partial_n)$ ($1\leq i\leq n$) is also constant
coefficient (i.e. $L\in\PD_0$).}

Let us suppose that $[p,L]=[P_i,L]=0$.  Then since $p\circ L=L\circ p$
and $p\circ\partial_i=\partial_i\circ p$ it follows that
\begin{eqnarray*}
0 = [P_i,L] &=& \partial_i\circ p\circ L-L\circ\partial_i\circ p\\
 &=& \partial_i\circ L\circ p -L\circ p\circ \partial_i\\
 &=& [L\circ p,\partial_i]
\end{eqnarray*}

Then we note that letting $L\circ p=\sum f_{\alpha}(x_1,\ldots,x_n)
\partial_1^{\alpha_1}\cdots \partial_n^{\alpha_n}$ be a series
representation for any microdifferential operator then
$$
[\partial_i,L]=\sum
f'_{\alpha}(x_1,\ldots,x_n)\partial_1^{\alpha_1}\cdots
\partial_n^{\alpha_n}
$$
where prime denotes differentiation with respect to $x_i$.  Hence, if
$L\circ p$ commutes with each $\partial_i$ then $L\circ p$ has constant
coefficients.  However, this provides a linear relation between
any coefficient of $L$ and certain higher ones.  Since the
coefficients of $L$ are bounded factorially with the order, this is
only possible if the coefficients are all constant.
\end{lemma}

\subsection{A Maximal Commutative Ring from Darboux Transformation}

Here we will
consider a special class of commutative rings of differential
operators for which we are able to demonstrate maximality using the
results of the previous subsection.
Suppose that the constant coefficient operator
$p\in\Dd_0$ factors as 
\begin{equation}
p(\partial_1,\ldots,\partial_n)=L\circ K\qquad
L,K\in\Dd. \label{eqn:factor}
\end{equation}
(See \bibref[BK] and \bibref[jmaa] for a discussion of some methods
for achieving such factorizations.)
Then the
method of Darboux transformation commonly used in the study of
integrable systems \bibref[Andrianov,Darboux1,Darboux2,Matveev] is to consider the (more complicated) operator
$$P:=K\circ L=K\circ p\circ K^{-1}$$ which shares many features with $p$
since the two operators are conjugate.  For instance, one may try to
conjugate other constant coefficient operators by $K$ to produce
operators that commute with $P$.  In fact, given any constant
coefficient operator $r(\partial_1,\ldots,\partial_n)\in\Dd_0$ it
follows that $[K\circ r\circ K^{-1},P]=0$ .  However, although $K\circ
r\circ K^{-1}\in\PD$ there is no reason to expect that it is in $\Dd$.
The content of the next theorem is the statment that the ring of all
\textit{differential} operators commuting with $P$ which are produced in this
way is a maximal commutative ring.

\note[Notation] Let $K\in\Dd$ be a differential operator and define
$$
R_0(K):=\{r\in\Dd_0\ |\ K\circ r(\partial_1,\ldots,\partial_n)\in \Dd\circ
K\}\subset \Dd_0
$$
to be the subring of elements $r\in\Dd_0$ such that $K\circ r$ has $K$ as
a right factor.  Then the ring
\begin{eqnarray*}
R(K)&:=&\left(K\circ \Dd_0\circ K^{-1}\right)\cap \Dd\\
 &=& K\circ R_0(K)\circ K^{-1}
\end{eqnarray*}
is a commutative subring of $\Dd$.  In general, it will be the case
that $R(K)=\C$ is trivial, but if $K$ is chosen to be a non-constant
operator satisfying \eqref{eqn:factor} then $R(K)$ will contain
differential operators.

\begin{theorem}[thm:maximal]{Let $K\in\Dd$ be a differential operator which is the
right factor of some constant coefficient operator $p=L\circ K\in\Dd_0$.
Then if $R'$ is a commutative ring such that
$$
R(K)\subset R'\subset \Dd
$$
it follows that $R(K)=R'$.  In other words, $R(K)$ is a maximal
commutative subring of $\Dd$.}

Let us suppose that $Q\in\Dd$ commutes with every element of $R(K)$.  We
must show that $Q$ is already in $R(K)$.
Note that since $p=L\circ K$ one automatically has that
$P_i:=\partial_i\circ p$ is in $R_0(K)$.  Thus, we know that $$[K\circ
P_i\circ K^{-1},Q]=0
\qquad\textup{for}\qquad1\leq i\leq n$$
and hence conjugation in the ring of pseudo-differential operators gives
$$
[P_i,K^{-1}\circ Q\circ K]=0\qquad 1\leq i\leq n.
$$
By Lemma~\ref{lem:only-const}, this implies that 
$$
K^{-1}\circ Q\circ K\in \PD_0
$$
and so it is clear that
$Q$ is of the form $K\circ q(\partial_1,\ldots,\partial_n)\circ
K^{-1}$ for some constant coefficient pseudo-differential operator $q\in\PD_0$.

However, we may moreover note that $q\in\Dd_0\subset \PD_0$ is a
constant coefficient \textit{differential} operator.  To show this we
introduce a normalized common eigenfunction
$$\psi(x_1,\ldots,x_n,z_1,\ldots,z_n)=\frac{1}{g(z_1,\ldots,z_n)}K\left[e^{x_1z_1+\cdots
+x_nz_n}\right] $$ where the polynomial $g$ is to be defined below.
Note that regardless of the choice of $g$, we have by construction
that $$ K\circ L[\psi]=p(z_1,\ldots,z_n)\psi\qquad
Q[\psi]=q(z_1,\ldots,z_n)\psi.  $$

One can write the function $K[\exp(\sum x_iz_i)]$ in the form
$$
K\left[e^{x_1z_1+\cdots
+x_nz_n}\right]=\left(\frac{\sum_{\alpha=1}^{N}
\rho_{\alpha}(z_1,\ldots,z_n) \sigma_{\alpha}(x_1,\ldots,x_n)}{\sigma_0(x_1,\ldots,x_n)}\right) e^{x_1z_1+\cdots
+x_nz_n}
$$
where $\rho_{\alpha}$  are all non-zero polynomials in $z_1,\ldots,z_n$
and $\sigma_{\alpha}$  are \textit{distinct, non-zero monomials} in
$x_1,\ldots,x_n$. We choose $g(z_1,\ldots,z_n)$ to be the highest
common factor of the polynomials $\rho_{\alpha}(z_1,\ldots,z_n)$.  We
have thus constructed $\psi$ so the product $f(z_1,\ldots,z_n)\psi$ is
holomorphic in each $z_i$ for $f\in \C((z_1,\ldots,z_n)$ if and only
if $f\in\C[z_1,\ldots,z_n]$ is actually a polynomial.

Then notice that for any $M\in\Dd$ one still has that 
$M[\psi]$ is holomorphic in $z_i$.  In particular, $M[\psi]$ is always
a polynomial in $\C(x_1,\ldots,x_n)[z_1,\ldots,z_n]$ multiplied by the
exponential function $\exp\sum x_iz_i$.  Putting this all together, since we have
already seen that 
$Q[\psi]=q(z_1,\ldots,z_n)\psi
$,
we conclude that $q(\partial_1,\ldots,\partial_n)\in\Dd_0$ is a
constant coefficient differential operator.

Finally, if $q\in \Dd_0$ has the property that $K\circ q\circ
K^{-1}\in \Dd$ then this implies that $q\in R_0(K)$ and hence that $Q\in
R(K)$ in the first place.  So no operators outside of $R(K)$ have the
property that they commute with every element of $R(K)$.
\end{theorem}

\subsection{An Explicit Example}

It is still not clear from the theorem above whether the maximal
commutative subrings of $\Dd$ constructed by Darboux transformation
require an infinite number of generators over $\C$.  Certainly in
trivial cases (e.g. $K\in\Dd_0$) the resulting ring may require only a
finite number of generators.   But it would be nice to prove that this
is always the case or alternatively to observe at
least one example which does not.  By considering a particular example
in detail here, we achieve the latter.

\subsubsection{A Subring of $\C[x,y]$}

\note[Notation] Let $\R\subset\C[x,y]$ ($\lambda\in\C$) be the subset 
$$
\R=\{q(x,y)\in\C[x,y]\ |\
q_x(z,\frac{\lambda}{z})=q_y(z,\frac{\lambda}{z})=q_{xy}(z,\frac{\lambda}{z})\equiv
0\}.  $$ In other words, $\R$ is the set of
polynomials $q\in\C[x,y]$ such that $q_x$, $q_y$ and
$q_{xy}$ all have a factor of $xy-\lambda$.  Note that $\C\subset\R$
and more importantly that if $q_1,q_2\in\R$ are two such polynomials
then $q_1+q_2\in\R$ and $q_1q_2\in\R$.  This obviously gives us that
\statelemma[lem:the-ring]{$\R$ is a proper subring of $\C[x,y]$ containing
$\C$ as well as every polynomial of the form $\rho(x,y)(xy-\lambda)^3$
for $\rho\in\C[x,y]$.}

It will be shown below that for a particular choice of $K\in\Dd$ the
maximal commutative subring $R(K)\subset\Dd$ is isomorphic to $\R$.
Therefore it is interesting to note that this ring requires an
infinite number of generators over $\C$.

\begin{lemma}[lem:inf-num]{The ring $\R$
has the form  $\C[\omega_1 (xy-\lambda)^3,\omega_2(xy-\lambda)^3,\ldots]$ where
$\{\omega_i\}$ is any basis of $\C[x,y]$ as a vector space.  In
particular, $\R$ is not finitely generated.}

We must show that a polynomial $q\in \C[x,y]$ is in $\R$ if and only
if it is of the form
$g(x,y)(xy-\lambda)^3+c$ for some $g\in\C[x,y]$ and $c\in\C$.
Clearly, such a $q$ is an element of $\R$.  Alternatively, let us
suppose that $q\in\R$ and therefore $q_x=(xy-\lambda)r(x,y)$.  Then,
since
$$
q_{xy}=xr(x,y)+(xy-\lambda)r_y(x,y)
$$
also has a factor of $xy-\lambda$ one finds that $r$ has a factor of
$xy-\lambda$ and hence $q_x$ actually has a factor of
$(xy-\lambda)^2$.  (Similarly for $q_y$.)

Now we have that $q_x=(xy-\lambda)^2g(x,y)$ for some $g\in\C[x,y]$.  Then
integrating by parts with respect to $x$ one has
$$
q(x,y)=\frac{1}{3y}\left[(xy-\lambda)^3g(x,y)-\int (xy-\lambda)^3
g_x(x,y)\,dx\right].
$$
Continuing to integrate by parts (choosing always to integrate
$(xy-\lambda)^j$ so that one gets higher powers of $xy-\lambda$ and
higher derivatives of $g$) one 
gets a finite sum (since a high enough derivative of $g$ will
eventually vanish) of terms each having a factor of $(xy-\lambda)^3$,
plus a constant of integration at the end.

Now we note that $\R$ cannot be constructed by a finite number of
generators over $\C$.  We may, for instance, suppose that
$R=\C[\nu_1,\nu_2,\ldots]$ for some polynomials $\nu_i$ and
w.l.o.g. we may take $\nu_i$ to have no constant term
($\nu_i(0,0)=0$).  But then consider the polynomials
$x^i(xy-\lambda)^3$ which are elements of $\R$.  These cannot involve
any \textit{products} of the generators $\nu_i$ else they would have a
factor of $xy-\lambda$ to a higher degree.  Thus, they must be a
\textit{linear} combination of some generators $\nu_i$.  On the other hand,
 since the polynomials $x^i(xy-\lambda)^3$ are linearly independent for
different $i$'s, it follows that you need infinitely many generators
to construct $\R$.
\end{lemma}

\subsubsection{Isomorphism to a ring of Differential Operators}
\label{subsec:R(K)}

Let us use the notation of the preceding subsection to describe a maximal
commutative ring of differential operators.
For this example, we will be working in two dimensions only, so $n=2$.
The constant coefficient differential operator which we will factor is
$p(\partial_1,\partial_2)=(\partial_1\partial_2-\lambda)^3$
($\lambda\in\C$) which factors as $p=L\circ K$ with
$$K=x_1x_2(\partial_1\partial_2-\lambda)\circ\frac{1}{x_1x_2}.$$
and
\begin{eqnarray*}
L &=& {{\partial_{1}}^2} {{\partial_{2}}^2}
+ {\frac{1}{x_{1}}}  \partial_{1} {{\partial_{2}}^2}
 -{{x_{1}}^{-2}}  {{\partial_{2}}^2}
+ {\frac{1}{x_{2}}}  {{\partial_{1}}^2} \partial_{2}\\
&&+ {\frac{1 - 2 \lambda x_{1} x_{2}}{x_{1} x_{2}}} 
   \partial_{1} \partial_{2}
+ {\frac{-1 - \lambda x_{1} x_{2}}{{{x_{1}}^2} x_{2}}}  \partial_{2}
 -{{x_{2}}^{-2}}  {{\partial_{1}}^2}
+ {\frac{-1 - \lambda x_{1} x_{2}}{x_{1} {{x_{2}}^2}}}  \partial_{1}\\
&&
+ {{\lambda}^2} + {\frac{1}{{{x_{1}}^2} {{x_{2}}^2}}} +
    {\frac{\lambda}{x_{1} x_{2}}}
\end{eqnarray*}

\begin{lemma}[lem:kernel]{A constant coefficient operator
$q(\partial_1,\ldots,\partial_n)\in\Dd_0$ is an element of $R_0(K)$ 
if and only if the function
$$
\psi(x_1,x_2,z):=x_1x_2 e^{x_1z +x_2\frac{\lambda}{z}}
$$
is in the kernel of the operator $K\circ q$ for all values of $z\in\C$.
}

One direction is especially
simple.  If $K\circ q=Q\circ K$ then
$$
K\circ q[x_1x_2 e^{x_1z+x_2\frac{\lambda}{z}}]=Q\circ x_1x_2
(\partial_1\partial_2-\lambda)[e^{x_1z+x_2\frac{\lambda}{z}}]\equiv 0.
$$

Conversely, let us suppose that $K\circ q$ annihilates this function.  This
means that $M:=K\circ q\circ {x_1x_2}$ applied to
$exp(x_1z_1+x_2z_2)$ is zero for all $z_1z_2-\lambda=0$.  But note
that $M$ applied to this exponential results in a polynomial in $z_i$
with coefficients in $\C(x_1,x_2)$ multiplied by an exponential.  This
product vanishes on $z_1z_2-\lambda=0$ if and only if the polynomial has
a factor of $z_1z_2-\lambda$ which implies that 
$$
M=Q\circ x_1x_2 (\partial_1\partial_2-\lambda)
$$
for \textit{some} $Q\in\Dd$.
Multiplying this equation on the right by $\frac{1}{x_1x_2}$ on the right proves the lemma.
\end{lemma}

Using this lemma and the previous theorem, as well as the
bispectrality \bibref[DGr,G] of the constant coefficient operators, we
demonstrate an isomorphism between $R(K)$ and $\R$.
\begin{theorem}[thm:isom]{
The ring $R(K)$, known to be maximal commutative by the preceding
theorem, is isomorphic to the ring $\R$ (cf.~Lemma~\ref{lem:the-ring}). }

Using Lemma~\ref{lem:kernel} and Theorem~\ref{thm:maximal}, we know that 
$R(K)$ is isomorphic to the ring
$$
R_0(K)=\{q\in \Dd_0\ |\ K\circ q[x_1x_2 e^{x_1z+x_2\lambda{z}}]\equiv0\}.
$$
However, since $\partial_{z_i}:=\frac{\partial}{\partial z_i}$
commutes with differential operators in the variables $x_i$, this
property is equivalent to saying that
$$
\partial_{z_1}\partial_{z_2}\left[ K\circ q
[ e^{x_1z_2+x_2z_2}]\right]\equiv 0\qquad \forall z_1z_2-\lambda=0.
$$
This can be written as differential equations for $q$ by applying all
of these differential operators, clearing the denominator by
multiplying by a polynomial in $x_1,x_2$ and looking at the
coefficients of each monomial in $x_i$.  These will be differential
expressions for polynomials in $z_1$ and $z_2$ including $q$ which
must vanish on $z_1z_2-\lambda$.  For this to happen, it is necessary
and sufficient that $q_x$, $q_y$ and $q_{xy}$ all have
$z_1z_2-\lambda$ as a factor.
\end{theorem}

\section{Resultants of Commuting Differential Operators}

In this section we give a definition of resultants for partial
differential operators (cf.~\bibref[CarraFerro]) including ``spectral
parameters'' and their significance in the commutative case.  

\subsection{Definitions}

Fix $0<n\in\N$ and denote by $\Omega^d$ the ${n+d\choose n}$-component
vector 
$$
\Omega^d=(\omega_1^d,\omega_2^d,\ldots)
$$
where $\omega_i^d$ run over all monomial, monic differential operators
in the variables $x_1,\ldots,x_n$ of degree less than or equal to $d$.
In other words, 
$$\omega_i^d\in\{\partial_1^{\alpha_1}\cdots \partial_n^{\alpha_n}\ :\
\alpha_i\in\N,\ \sum \alpha_i\leq d\}
$$
By writing it as a vector, we are supposing that they have an
ordering.  Such an ordering is a choice, not determined canonically,
but the particular choice is not important to the following.
Then, for any differential operator $L$ of order $d$ or less, we
denote by $\v_d{(L)}$ the vector whose $i^{th}$ entry is the coefficient
of $\omega_i^d$ in $L$.  In particular, $L=\v_d{(L)}\cdot \Omega^d$.
Let $L_1,\ldots,L_{n+1}$ be differential operators in the variables
$x_1,\ldots, x_n$ having orders $l_1,\ldots,l_{n+1}$ respectively.

Let $N:=-n+\sum l_i$ and construct the matrix
$R_{\mu}=R_{\mu}(L_1,\ldots,L_{n+1})$ whose rows are
$\v_{N}{(\omega_j^{N-l_i}\circ ( L_i-\mu_i))}$ for all $1\leq i\leq
n+1$ and all $1\leq j\leq {n+N-l_i\choose n}$.  We call any maximal
minor determinant of $R_{\mu}$ a \textit{partial $\mu$-shifted
differential resultant}.  Note that each partial $\mu$-shifted
differential resultant is a polynomial in the variables $\mu_i$
($1\leq i\leq n+1$) with coefficients that may depend on $x_j$ ($1\leq
j\leq n$).  We define the \textit{$\mu$-shifted differential
resultant} of the operators $L_i$ to be the polynomial in the
variables $\mu_i$ which is the greatest common divisor of all of these
maximal minor determinants.

\note In the case $L_i\in\C[\partial_1,\ldots,\partial_n]$, this
definition is a special case of the polynomial resultant
\bibref[Macaulay] and that in one dimension with variable coefficients
it reproduces the differential resultant of ordinary differential
operators used to construct the spectral curve
\bibref[PrevWeil,PrevSpec].  The definition of the
\textit{differential resultant} of the operators $L_i$ given in
\bibref[CarraFerro] is, in our terminology, a particular
\textit{partial} $\mu$-shifted differential resultant of the operators
$L_i$ with all $\mu_i=0$.

\subsection{The Commutative Case}

As in the one dimensional case, we will here show that the
differential resultant provides a polynomial equation satisfied by the
operators $L_i$ in the case that they mutually commute.
The remainder of the section will then be comprised of examples and
counter-examples of what we would hope to have as a consequence.
First, following the approach used in \bibref[CarraFerro], we
demonstrate the following essential lemma:

\begin{lemma}[lem:dform]{Any partial $\mu$-shifted
differential resultant of the operators $L_1,\ldots,L_{n+1}$ can be
written as
\begin{equation}
\sum_{i=1}^{n+1} D_i\circ(L_i-\mu_i)\label{eqn:dform}
\end{equation}
for some partial differential operators $D_i$ with coefficients
depending on $\mu_j$ ($1\leq j\leq n+1$) and $x_k$ ($1\leq k\leq n$).}

Let $j$ be the integer $1\leq j\leq {n+N\choose n}$ such that
$\omega_j^N=1$ is the differential operator of order zero in the
vector $\Omega^N$.
Construct the matrix $M$ of size
${n+N\choose n}\times{n+N\choose n}$ which is the identity matrix
except for the fact that the $j^{th}$ column is replaced by the vector
$\Omega^N$.  Note that $\det M=1$.

Let $\tilde R$ be a maximal square minor of the matrix
$R_{\mu}(L_1,\ldots,L_{n+1})$.  Note that the elements in the $j^{th}$
column of the matrix $\tilde R\cdot M$ are all monic monomial
differential operators composed with the operators $L_i-\mu_i$.  Then,
expanding down this column while taking determinants, one finds
exactly something of the form \eqref{eqn:dform} with the coefficients
of $D_i$ coming from the other minor determinants of $\tilde R$.

On the other hand, it is an elementary fact of linear algebra that
$\det\tilde R=\det \tilde R\cdot M$ and so \eqref{eqn:dform} must
actually be equal to the order zero operator which is the partial
$\mu$-shifted differential resultant of the operators $L_i$.
\end{lemma}

Now suppose that the operators $L_i$ ($1\leq i\leq n+1$) mutually
commute.  By definition, any partial $\mu$-shifted differential
resultant of these operators is a polynomial
in the variables $\mu_i$ with coefficients possibly depending on the
variables $x_j$.  As a consequence of Lemma~\ref{lem:dform} we then
find that the operators $L_i$ \textit{satisfy} this polynomial.

\begin{theorem}[thm:satisfy]{Let  $p(\mu_1,\ldots,m_{n+1})$ be any
$\mu$-shifted differential resultant of the mutually commuting
operators $L_i$, then 
$$
p(L_1,\ldots,L_{n+1})=0.
$$}
Only commutativity of the $\mu_i$'s with the $L_j$'s is required to
rewrite $p(\mu_1,\ldots,\mu_{n+1})$ in the form \eqref{eqn:dform}.
So, since $[L_i,L_j]=0$ we can write $p(L_1,\ldots,L_{n+1})$ by
substituting $L_i$ for $\mu_i$ in \eqref{eqn:dform}.  This, however,
is clearly zero since every term has a factor of $L_i-\mu_i$ for some
$i$.  
\end{theorem}

In the one dimensional case, we can moreover say that the
$\mu$-shifted differential resultant is a polynomial in $\mu_1$ and
$\mu_2$ with \textit{constant} coefficients or a multiple of such a
polynomial by a function of $x_1$.  Here, the results proved thus far
leave open the possibility that the differential resultant will only
produce polynomial equations satisfied by the operators with
\textit{explicit} dependence on the variables $x_j$.  We were not able
to produce any such examples or exclude the possibility. 

\subsection{The Zero Possibility}

As stated in the introduction, given two commuting ordinary
differential operators $L_1$ and $L_2$ the determinant of $R_{\mu}(L_1,L_2)$
(which happens to always be square in the case $n=1$) is a
\textit{non-zero} polynomial in $\mu_1$ and $\mu_2$ which is satisfied
by the operators.  Here we will see that the differential resultant
does not always give such useful information in the higher dimensional
case. 

Consider the case $n=2$ and 
$$
L_1=\partial_1^2-\partial_2^2-1\qquad 
L_2=\partial_1\circ L_1\qquad L_3=\partial_2\circ L_1.
$$
Note that these operators satisfy the equation
$L_2^2-l_3^2-L_1-L_1^6=0$ and so one might hope, given
Theorem~\ref{thm:satisfy}, that the differential resultant of these
operators is $\mu_2^2-\mu_3^2-\mu_1-\mu_1^6$ (or at least is a
non-zero multiple of this).

\begin{lemma}{The differential resultant of the operators $L_i$ is
the zero polynomial in the variables $\mu_i$ ($1\leq i\leq 3$).} One
could, of course, merely compute the resultant according to the
definition.  However, there is a more direct and informative way to
observe this fact.  Since these operators are constant coefficient,
the problem reduces to a problem of polynomial resultants.  In
particular, the resultant is the same as the resultant of the
homogeneous polynomials $$
p_1(x_1,x_2,x_3)=x_1^2-x_2^2-(\mu_1+1)x_3^2\qquad
p_2(x_1,x_2,x_3)=x_1^3-x_1x_2^2-x_1x_3^2-\mu_2x_3^3$$
$$
p_3(x_1,x_2,x_3)=x_1^2x_2-x_2^3-x_2x_3^2-\mu_3x_3^3
$$
However, it is well known \bibref[Macaulay] that this resultant will
be zero iff these polynomials have a common zero in projective space.
Although it is true that no ``finite'' point (with $x_3\not=0$) is a
common solution to these polynomials for all values of $\mu_i$, there
are solutions at infinity.  In particular,
note that the point $(1,-1,0)$ satisfies all three polynomials
regardless of the values of $\mu_i$. 
\end{lemma}

It is interesting to note the geometry behind this situation.  This
problem of having a zero resultant never arises in the
one dimensional case essentially because only one point is being added
at infinity and that point is never a solution of the homogeneous
polynomial.  Whereas, in higher dimensions, there is ``room'' at
infinity for many solutions.

Note that the same problem can also occur in a non-constant case (and
so not simply an example of a polynomial resultant).  In particular,
the differential resultant of any three operators from the ring $R(K)$
described in Subsection~\ref{subsec:R(K)} will be zero regardless of
the values of the variables $\mu_i$.  The mundane explanation of this
fact here is merely that the $N^{th}$ powers of $\partial_1$ and
$\partial_2$ never appear in $\omega_j^{N-l_i}\circ ( L_i-\mu_i)$ and
so there are columns of the resultant matrix with all zero entries.

\subsection{Positive Results}

A more encouraging example is to consider the operators
$$
L_1=\partial_1^2-\partial_2^2
\qquad
L_2=x_2\partial_1+x_1\partial_2
\qquad
L_3=L_1\circ L_2 - \gamma L_1\qquad \gamma\in\C.$$
It is a non-obvious fact that $[L_1,L_2]=0$, but given this (which is
easily checked) it is clear that $L_3$ also commute and that the three
together satisfy a polynomial equation $p(L_1,L_2,L_3)=0$ with
$$
p(\mu_1,\mu_2,\mu_3)=\mu_3-\mu_1\mu_2+\gamma \mu_1.
$$
Then the differential resultant (which can be most easily computed not
by finding all maximal minor determinants but by the formula $D/A$
where $D$ and $A$ are the minor determinants specified in
\bibref[Macaulay]) is exactly $p^3(\mu_1,\mu_2,\mu_3)$.   This is very
nearly what we would want (although there is presently no theory to
explain the exponent ``3'' which arises).

It is intriguing and surprising that the resultant is independent of
the variables $x_1$ and $x_2$ in this case.  In the one dimensional
case, the resultant of two monic differential operators is independent
of $x$ if and only if the operators commute.  Here, the situation
involves one operator, $L_1$, with constant leading coefficients and
others that are not, which cannot happen in the one dimensional case.

\section{Mathematica Code}

The following code (and assistance using it) is available by writing
to the first author ({\tt kasman@math.cofc.edu}).  It is useful for
performing many calculations with differential operators, and here we
will only give a very brief description of how to use it.
First set the value of the variable
{\tt dimen} equal to the number of variables you will be working with
and then input the file containing the code.  For instance:

{\small\begin{verbatim}

In[1]:= dimen=3

Out[1]= 3

In[2]:= <<pdo-ak.m
String dimen already defined..using present value 3
READY: Partial Differential Operators in 3 dimensions
\end{verbatim}}

Then, a differential operator can be entered using {\tt
DX[1], DX[2],...} as the elementary differential operators and {\tt
X[1], X[2],...} as the corresponding variables.  For example:
{\small\begin{verbatim}


In[3]:= L=X[1] DX[2]+X[2] DX[1]

Out[3]= DX[2] X[1] + DX[1] X[2]

In[4]:= Q=DX[1]^2+DX[2]^2

             2        2
Out[4]= DX[1]  + DX[2]

\end{verbatim}}

Note that when entering or reading a differential operator in this
notation, it is always assumed that all differentiation has been taken
and functions are on the ``left'', even if it may not be written this
way.  (That is, {\tt X[1] DX[1]=DX[1] X[1]}.)  The non-commutativity
is only apparent when multiplying two differential operators.
Multiplication of differential operators is achieved using the command
{\tt pdomult}:
{\small\begin{verbatim}
In[6]:= pdomult[DX[1],X[1]]

Out[6]= DX[1] X[1] + 1

In[7]:= pdomult[Q,L]
                            2
18:10 > multiplying by DX[1]
                            2
18:10 > multiplying by DX[2]
Simplifying...
.
.
.
18:10 > Simplifying {3, 0, 0} term.
.
18:10 > Simplifying {1, 1, 0} term.
18:10 > Simplifying {2, 1, 0} term.
.
.
18:10 > Simplifying {1, 2, 0} term.
.
.
18:10 > Simplifying {0, 3, 0} term.
.
.
.

                             2                   3             3
Out[7]= 4 DX[1] DX[2] + DX[1]  DX[2] X[1] + DX[2]  X[1] + DX[1]  X[2] + 
 
                2
>    DX[1] DX[2]  X[2]

\end{verbatim}}
Similarly, if you want to \textit{apply} the operator {\tt L} to the
function {\tt f} you simply say {\tt pdoapply[L,f]}.

Finally, given a list of $N$ operators in $N-1$ variables, the matrix
whose minor determinants give the differential resultant can be
computed as {\tt DiffResult[{L1,L2,...,LN}]}.  (Note that this command
automatically subtracts the indeterminate {\tt mu[i]} from the
$i^{th}$ operator so that the differential resultant is a polynomial
in these variables.  The command {\tt
TakeRandomDeterminants[mat,n]} may come in handy as well, since it
takes the determinant of {\tt n} randomly chosen maximal minors of the
matrix {\tt mat}.

Here is the code:
{\small\begin{verbatim}
(* PDO-AK.M  Mathematica Code for Differential Operators
by Alex Kasman, College of Charleston 
kasman@math.cofc.edu *)

(* Set dimen to be the number of dimensions and then input this
file.  It will be 2 by default *)

(* Use X[1] ... X[dimen] as variables and DX[1] ... DX[dimen] as the
corresponding elementary differential operators. *)

(* Multiply two operators with pdomult[L,Q] *)

(* Conjugate one operator by another by pdoconj[L,Q], which yields the
operator Q L Q^1 if this is a differential operator. *)

(* DiffResult[{L1,....,LN}] gives the matrix whose minor determinants
have a gcd which is the differential resultant of L1-mu[1], L2-mu[2],
etc. *)

If[StringMatchQ[ToString[dimen],"dimen"],
dimen=2;
Print["String dimen undefined...setting to ",dimen," by default..."],
Print["String dimen already defined..using present value ",dimen]]


(* Call the variables X[1], X[2],...X[dimen] and the differential operators
DX[1],DX[2],...DX[dimen] *)

zerovect=Table[0,{i,1,dimen}]

DXpower[vect_]:=Module[{i},Product[DX[i]^(vect[[i]]),{i,1,Length[vect]}]]

DXv[alpha_]:=Module[{i},Table[DX[i],{i,1,Length[alpha]}] . alpha]

pdocoef[L_,zerovect]:=L/. DX[i_]->0

pdocoef[L_,vect_]:=Coefficient[Collect[((Expand[L]
/.DXpower[vect]->SPACE)/.DX[i_]->0),SPACE],SPACE]

pdotermmult[coef_,vect_,M_]:=Module[{i,sofar},
verbose["multiplying by ",DXpower[vect]];
For[i=1;sofar=M,i<Length[vect]+1,i=i+1,
sofar=pdopowmult[i,vect[[i]],sofar]];
coef sofar]

pdotermmult[0,vect_,M_]:=0

pdopowmult2[i_,1,M_]:=D[M,X[i]]+M DX[i]
pdopowmult[i_,0,M_]:=M

pdopowmult[i_,n_,M_]:=(pdopowmult2[i,n,M])
pdopowmult2[i_,n_,M_]:=pdopowmult2[i,n-1,pdopowmult2[i,1,M]]

maketablefor[M_]:=Module[{max,i,j,k,sofar},
For[i=1,i<=dimen+1,i=i+1,max[i]=Exponent[Collect[M,DX[i]],DX[i]]];
For[i=1;sofar=Table[j[k],{k,1,dimen}],i<dimen+1,i=i+1,
sofar=Table[sofar,{j[i],0,max[i]}]];
Flatten[sofar,dimen-1]]

maketablefor[M_,K_]:=Module[{max,i,j,k,sofar},
For[i=1,i<=dimen+1,i=i+1,max[i]=Exponent[Collect[M,DX[i]],DX[i]]
-Exponent[Collect[K,DX[i]],DX[i]]];
For[i=1;sofar=Table[j[k],{k,1,dimen}],i<dimen+1,i=i+1,
sofar=Table[sofar,{j[i],0,max[i]}]];
Flatten[sofar,dimen-1]]

pdoexp[ll_,n_]:=Module[{outL,j},(verbose["Raising operator to power
",n]; For[j=0;outL=1,j<n,j=j+1,verbose["pdoexp: power
",j];outL=pdomult[ll,outL]];outL)]

makeshortlist[L_]:=Module[{i,inlist,outlist}, verbose["making short
list"]; outlist={}; inlist=maketablefor[L];
For[i=1,i<Length[inlist]+1,i=i+1,
If[MatchQ[ToString[pdocoef[L,inlist[[i]]]],"0"],verbose["not in
list",inlist[[i]]],outlist=Union[outlist,{inlist[[i]]}];
verbose["makeshortlist: includes ",inlist[[i]]]]]; outlist]

pdomult[L_,M_]:=Module[{i,j,LM,liszt},
liszt=maketablefor[L];
For[i=1;LM=0,i<=Length[liszt],i=i+1,
LM=LM+pdotermmult[pdocoef[L,liszt[[i]]],liszt[[i]],M]];
pdosimp[LM]]

pdotermapply[coef_,vect_,M_]:=Module[{i,sofar},
verbose["Applying ",DXpower[vect]];
For[i=1;sofar=M,i<Length[vect]+1,i=i+1,
sofar=pdopowapply[i,vect[[i]],sofar]];
coef sofar]

pdotermapply[0,vect_,M_]:=0

pdopowapply[i_,n_,M_]:=(D[M,{X[i],n}])

pdoapply[L_,M_]:=Module[{i,j,LM,liszt},
liszt=maketablefor[L];
For[i=1;LM=0,i<=Length[liszt],i=i+1,
LM=LM+pdotermapply[pdocoef[L,liszt[[i]]],liszt[[i]],M]];
Simplify[LM]]

pdosimp[0]:=0

pdosimp[L_]:=Module[{liszt,i,sofar}, Print["Simplifying..."];
liszt=maketablefor[L]; For[i=1;sofar=0,i<Length[liszt]+1,i=i+1,
If[MatchQ[ToString[pdocoef[L,liszt[[i]]]],"0"],Print["."],
verbose["Simplifying ",liszt[[i]]," term."];
sofar=sofar+Simplify[pdocoef[L,liszt[[i]]]] DXpower[liszt[[i]]]]];
sofar]

pdodisplay[L_]:=Module[{liszt,i},
liszt=maketablefor[L];
For[i=Length[liszt],i>0,i=i-1,
pdotermdisplay[pdocoef[L,liszt[[i]]],liszt[[i]]]]]

pdotermdisplay[0,vect_]:=(bleh=0)

pdotermdisplay[coef_,vect_]:=Print["+ (",Simplify[coef],") ",DXpower[vect]]

pdoconj[P_,K_]:=Module[{L,lco,LK,KP,i,j,liszt,deg,bleh,slv,holdit},
liszt=maketablefor[P]; L=Sum[lco[liszt[[i]]]
DXpower[liszt[[i]]],{i,1,Length[liszt]}]; deg=pdodeg[P];
For[i=1,i<Length[liszt]+1,i=i+1,
If[Sum[liszt[[i]][[j]],{j,1,dimen}]>deg-1, bleh[liszt[[i]]]=0;
L=L/.lco[liszt[[i]]]->pdocoef[P,liszt[[i]]], bleh[liszt[[i]]]=100]];
Print["finding KP"]; KP=pdomult[K,P]; Print["finding LK"];
LK=pdomult[L,K]; brahms=maketablefor[LK];
For[i=Length[liszt],i>0,i=i-1, verbose["solving for
lco[",liszt[[i]],"]"]; If[bleh[liszt[[i]]]==0,verbose["...Already
set..."], For[j=Length[brahms],j>0,j=j-1,
holdit=pdocoef[LK,brahms[[j]]];
If[StringMatchQ[ToString[D[holdit,lco[liszt[[i]]]]],"0"],Print["."],
slv=Solve[pdocoef[LK,
brahms[[j]]]==pdocoef[KP,brahms[[j]]],lco[liszt[[i]]]][[1]];
Print["replacement ->",slv]; j=-5; L=L/.slv; LK=LK/.slv]]]];
(*Print["Okay if this is zero ",pdosimp[LK-KP]];*) L=pdosimp[L]]

pdodeg[L_]:=Module[{transL,i,j,deg,const},
transL=L/.DX[i_]->const[i] SPACE;
deg=Exponent[Collect[transL,SPACE],SPACE];
deg]

pdoTeX[L_]:=Module[{liszt,i},
liszt=maketablefor[L];
For[i=Length[liszt],i>0,i=i-1,
pdotermTeX[pdocoef[L,liszt[[i]]],liszt[[i]]]]]

pdotermTeX[0,vect_]:=(bleh=0)
pdotermTeX[1,vect_]:=Print["+ ",TeXForm[DXpower[vect]]]
pdotermTeX[coef_,zerovect]:=Print["+ ",TeXForm[coef]]

pdotermTeX[coef_,vect_]:=Print["+ ",TeXForm[Together[coef]]," **
",TeXForm[DXpower[vect]]]

Print["READY: Partial Differential Operators in ",dimen," dimensions"]

genminusone[deg_,liszt_]:=Module[{list2,i,j,n},
list2={};
For[i=1,i<=Length[liszt],i=i+1,
n=Sum[liszt[[i]][[j]],{j,1,Length[liszt[[i]]]}];
list2=Flatten[{list2,Table[Flatten[{liszt[[i]],j},1],
{j,0,deg-n}]},1]];
list2]

degcomplete[deg_,entry_]:=Module[{i},
Flatten[{entry,deg-Sum[entry[[i]],
{i,1,Length[entry]}]},1]]

AllMonomialsDegExactly[deg_]:=Module[{list2,list3,i,j,n},
If[dimen==1,DX[1]^deg,
list2=Table[{i},{i,0,deg}];
For[i=1,i<dimen-1,i=i+1,list2=genminusone[deg,list2]];
Table[degcomplete[deg,list2[[i]]],{i,1,Length[list2]}]]]

AllMonomialsDegAtMost[deg_]:=Module[{list2,i,j,n},
list2=Table[{i},{i,0,deg}];
For[i=1,i<dimen,i=i+1,list2=genminusone[deg,list2]];
list2];

multbynextdxupto[deg_,liszt_]:=Module[{list2,list3,i,j,k,n,n2,n3,n4,op,rem},
list2={}; For[i=1,i<=Length[liszt],i=i+1, op=liszt[[i]][[1]];
list3=liszt[[i]][[2]]; n=Length[list3]; k=Sum[list3[[j]],{j,1,n}];
rem=op; list2=Flatten[{list2,{{rem,Flatten[{list3,{0}},1]}}},1];
For[j=1,j<=deg-k,j=j+1, rem=Simplify[D[rem,X[n+1]]+rem*DX[n+1]];
list2=Flatten[{list2,{{rem,Flatten[{list3,{j}},1]}}},1]]]; list2]

MultByAllMonomialsUpto[deg_,L_]:=Module[{i,j,k,list2,rem},
rem={{L,{0}}}; verbose["MultByAllMonomials j=",1,"/",dimen];
For[i=1,i<deg+1,i=i+1, rem=Flatten[ {rem,
{{D[rem[[Length[rem]]][[1]],X[1]]+rem[[Length[rem]]][[1]]
DX[1],rem[[Length[rem]]][[2]]+1}}}, 1]]; For[j=2,j<=dimen,j=j+1,
verbose["MultByAllMonomials j=",j,"/",dimen];
rem=multbynextdxupto[deg,rem]];
Table[Expand[rem[[k]][[1]]],{k,1,Length[rem]}]]

pwr[liszt_]:=Product[DX[i]^liszt[[i]],{i,1,Length[liszt]}]

DiffResult[liszt_]:=Module[{i,j,k,maxdeg,n,n2,n3,mm,list2,list3,list4,rem},
n=Length[liszt];
maxdeg=Sum[pdodeg[liszt[[i]]],{i,1,n}]-dimen;
mm={};
list4={};
list2=AllMonomialsDegAtMost[maxdeg];
n2=Length[list2];
For[i=1,i<=n,i=i+1,
verbose["Diff result i=",i,"/",n];
list3=MultByAllMonomialsUpto[maxdeg-pdodeg[liszt[[i]]],liszt[[i]]-mu[i]];
mm=Flatten[{mm,list3},1]];
mm=makeamatrix[mm,list2];
verbose["The Matrix is:"];
If[Length[mm]+Length[mm[[1]]]>30,
verbose["TOO BIG TO PRINT"],
verbose[MatrixForm[mm]]];
(*TakeAllSquareMinors[mm,list4]*)
mm]

makeamatrix[liszta_,list2_]:=Module[{liszt,i,j,k,n,list3,list4,bleh},
wid=Length[list2]; list4=liszta/.DX[i_]->0; liszt=liszta-list4;
list3=Table[DXpower[list2[[i]]]->Table[Which[k==i,1,True,0],
{k,1,wid}],{i,Length[list2],2,-1}];
Simplify[(list4*bleh/.bleh->Table[Which[k==1,1,True,0],
{k,1,wid}])+(liszt/.list3)]]

TakeAllSquareMinors[mm_]:=Module[{liszt,
list2,i,j,k,n,ll,ww,deter,rem,minlist}, n=0; ll=Length[mm];
ww=Length[mm[[1]]]; If[ww==ll,{Det[mm]}, If[ww>ll,verbose["(Note:
Transposing the matrix for simplicity)"];
TakeAllSquareMinors[Transpose[mm]], For[k=0,k<=ll-ww,k=k+1,
verbose["Doing Young diagrams starting with ",k]; liszt={{k}};
For[i=1,i<ww,i=i+1, list2={}; For[j=1,j<=Length[liszt],j=j+1,
list2=Flatten[{list2,Table[Flatten[{liszt[[j]],
k},1],{k,0,liszt[[j]][[Length[liszt[[j]]]]]}]},1]]; liszt=list2];
verbose["Computing the corresponding determinants.  Listing nonzero
ones:"]; For[j=1,j<=Length[liszt],j=j+1,
deter=Simplify[Det[Table[mm[[ll-ww+i-liszt[[j]][[i]]]],{i,1,ww}]]];
If[StringMatchQ[ToString[deter],"0"],, n=n+1; rem[n]=deter; verbose["
******************* "]; verbose[" "]; verbose["Nonzero determinant
number ",n]; verbose[deter]]]]; Table[rem[j],{j,1,n}]]]]

time:=StringJoin[ToString[Date[][[4]]],":",ToString[Date[][[5]]]," > "]

verbose[f_]:=Print[time,f]
verbose[f_,g_]:=Print[time,f,g]
verbose[f_,g_,h_]:=Print[time,f,g,h]
verbose[f_,g_,h_,i_]:=Print[time,f,g,h,i]
verbose[f_,g_,h_,i_,j_]:=Print[time,f,g,h,i,j]
verbose[f_,g_,h_,i_,j_,k_]:=Print[time,f,g,h,i,j,k]

normp[a_,b_]:=Which[b<a,True,True,False] 

TakeRandomDeterminants[mat_,n_]:=Module[{ll,ww,randomyoung,i,bleh},
ll=Length[mat];
ww=Length[mat[[1]]];
verbose[ll,"x",ww," matrix"];
randomyoung:=Sort[Table[Floor[Random[]*(ll-ww+1)],{i,1,ww}],normp];
For[i=1,i<n+1,i=i+1,
verbose["Try number ",i];
lst=randomyoung;
verbose["Young is ",lst];
bleh=Simplify[Det[Table[mat[[ll-ww+i-lst[[i]]]],{i,1,ww}]]];
verbose["Determinant is -> ",bleh]]]
\end{verbatim}}

\end{document}